\documentclass[11pt]{article}

\usepackage[margin=1in]{geometry}
\usepackage[T1]{fontenc}
\usepackage[utf8]{inputenc}
\usepackage{amsmath,amssymb,amsthm}
\usepackage{graphicx}
\usepackage{subcaption}
\usepackage{xcolor}
\usepackage{hyperref}
\usepackage{xr-hyper}
\usepackage{comment}

\hypersetup{
hidelinks,
pdftitle={Magnitude-Based Features for Multispecies Spatial Data},
pdfauthor={Julia Sollberger, Joshua Bull, Sara Kalisnik, Bernadette Stolz}
}

\title{Magnitude-Based Features for Multispecies Spatial Data}

\author{
Julia Sollberger\thanks{Vrije Universiteit Amsterdam, \href{mailto:j.sollberger@vu.nl}{j.sollberger@vu.nl}}
\and
Joshua Bull\thanks{University of Oxford, \href{mailto:joshua.bull@maths.ox.ac.uk}{joshua.bull@maths.ox.ac.uk}}
\and
Sara Kališnik\thanks{Pennsylvania State University, \href{mailto:skalisnik@psu.edu}{skalisnik@psu.edu}. Sara Kališnik and Bernadette Stolz are joint senior authors.}
\and
Bernadette Stolz\thanks{Max Planck Institute of Biochemistry and Munich Center for Machine Learning, \href{mailto:stolz@biochem.mpg.de}{stolz@biochem.mpg.de}}
}

\date{\today}

\newtheorem{theorem}{Theorem}

\newtheorem{proposition}[theorem]{Proposition}

\theoremstyle{definition}
\newtheorem{definition}[theorem]{Definition}
\newtheorem{example}[theorem]{Example}

\theoremstyle{remark}
\newtheorem{remark}[theorem]{Remark}

\newcommand{\mum}{\ensuremath{\mu\mathrm{m}}}
\newcommand{\RR}{\mathbb{R}}

\begin{document}

\maketitle

\begin{abstract}
Multispecies spatial data arise in many applications where interactions between different entities are central to system behaviour, including biomedical imaging, geospatial analysis, and species ecology. Despite their importance, relatively few quantitative tools exist to capture such interactions. In this work, we propose magnitude-based features for the analysis of multispecies spatial data. Magnitude is a real-valued invariant of finite metric spaces that can be interpreted as an effective number of points, incorporating both spatial configuration and scale. We develop global and local magnitude feature vectors and demonstrate their utility on synthetic tumour microenvironment data, and in tissue microarray data from human colorectal cancer samples. Locally, the method identifies distinct neighbourhood types and reveals spatial heterogeneity; in the model, this includes radial patterns associated with different qualitative outcomes of the simulations, while in the real-world data it reflects the importance of tertiary lymphoid structure-like interactions between B and T cell populations. Globally, the approach recovers known classifications of long-term simulation outcomes across parameter regimes in synthetic data, and suggests important roles for CD4+ T cells and CD163+ macrophages in distinguishing patients with favourable Crohn's like reactions from unfavourable diffuse immune infiltration. Together, these results suggest that magnitude-based features provide a powerful and flexible tool for the analysis of multispecies spatial data.
\end{abstract}



\vspace{1cm}

\section{Introduction}

Recent advances in experimental technologies have greatly improved the ability to study biological systems in detail. Techniques such as single-cell and multi-omics methods can now measure many different molecular features of cells simultaneously, revealing previously hidden complexity and heterogeneity. Complementing these approaches, spatial biology captures locations of cells and molecules within tissues, preserving spatial context. Spatial biology encompasses a range of techniques that capture spatial information across a variety of molecular and cellular modalities in biological samples, enabling the study of how the spatial organisation of cells and molecules relates to structure and function within tissues \cite{carstens2024spatial, jacquemet2020cell, karimi2024method,lewis2021spatial, moffitt2022emerging, park2022spatial, vandereyken2023methods, xiaowei2021method}.

These advances generate increasingly rich datasets that require computational methods beyond summary statistics to fully capture their structure and complexity. A key challenge is to quantify not only the abundance of cells or cell types, but also their spatial organisation. This problem has roots in ecology, where Solow and Polasky introduced the notion of an `effective number of species' as a way to quantify biodiversity in a decision-making framework for conservation~\cite{Solow1994}. Their approach incorporated both the number of species and pairwise dissimilarities, formalised through three principles: diversity should increase with the number of species (monotonicity), should not increase when duplicating an existing species (twinning), and should increase as species become more dissimilar (monotonicity in distance). The resulting quantity coincides with a concept now called \emph{magnitude}.

In 2011, Tom Leinster reintroduced magnitude in a mathematical context by generalizing the Euler characteristic from categories to metric spaces via enriched category theory~\cite{L13}. Since then, magnitude has emerged as a geometric invariant that, under suitable conditions, encodes fundamental features such as dimension, perimeter, area, and volume~\cite{L13, LW13, Meckes2015, Willerton2014}. Beyond its origins in biodiversity~\cite{LC12, Solow1994}, magnitude has found applications in areas such as machine learning~\cite{magnitude-ML, torkamani2026magnitudedistancegeometricmeasure} and has also been studied in connection with persistent homology~\cite{GH21, OMALLEY2023107396, O18} within topological data analysis (TDA). In these settings, efficient computation of magnitude is essential, and can be achieved in practice using methods such as Cholesky decomposition~\cite{magnitude-ML}.

Its mathematical properties suggest that magnitude is well suited for quantifying spatial heterogeneity in complex biological systems such as the tumour microenvironment (TME). The TME refers to the local environment surrounding a tumour, including interacting components such as tumour, immune, and stromal cells, blood and lymphatic vessels, molecular signals, and extracellular matrix. Its architecture, cellular composition, and local interactions play a critical role in disease progression and response to treatment \cite{hanahan2012accessories,spatialYuan2016}. Tumours and their microenvironments are highly heterogeneous, both across patients and within individual tumours~\cite{bareham2024defining,lewis2021spatial,spatialYuan2016}. For example, spatial organisation at the tumour core differs markedly from that at the invasive edge and can reflect distinct transition states \cite{bareham2024defining}. Moreover, these environments evolve dynamically over time during disease progression and treatment.

In analogy to ecology, spatial neighbourhoods in the TME are often described as specialised niches, motivating the transfer of concepts and methods from ecological theory  \cite{spatialYuan2016}.  Importantly, identical proportions of cell types within a niche do not necessarily correspond to functionally similar environments, highlighting the need for principled measures of spatial organisation \cite{spatialYuan2016}. Spatial biology now enables the study of the TME across multiple modalities at unprecedented resolution. The spatial arrangement of cells governs interactions and signalling processes, making the quantification of spatial patterns and neighbourhood structure essential for understanding function and developing biomarkers \cite{bareham2024defining, carstens2024spatial, lewis2021spatial, moffitt2022emerging, park2022spatial, vandereyken2023methods, spatialYuan2016, bull_muspan_2025}. While computational approaches such as spatial statistics and machine learning exist, defining robust and interpretable metrics for spatial heterogeneity beyond summary statistics at a fixed scale remains a central challenge \cite{bareham2024defining, lewis2021spatial, bull_muspan_2025}. Recently, new approaches in TDA, referred to as \emph{relational} or \emph{chromatic} TDA, have been introduced \cite{di_Montesano_2025, natarajan2026topology, Stolz2024}. However, these methods are often computationally expensive.

Here, we apply magnitude at both local and global scales to two datasets. Primarily, we consider simulated data from a mathematical model of the TME \cite{Bull2023}. Our results show that magnitude recovers known classifications of parameter regimes associated with distinct long-term behaviours. Locally, it identifies distinct neighbourhood types and reveals spatial heterogeneity. We further validate this approach by calculating magnitude at both local and global scales on a dataset of 140 tissue microarray samples from human colorectal cancer samples imaged using CODEX spatial proteomics~\cite{schurch_coordinated_2020}. Magnitude recapitulates key cell types known to be positively correlated with survival, while also identifying additional combinations of cell types which correlate with poor prognosis.

The main contributions of this work are as follows:
\begin{itemize}
\item We introduce magnitude-based feature vectors for multispecies spatial data. A key property underlying this construction is that magnitude does not satisfy the inclusion–exclusion principle, which enables it to capture spatial interaction effects.
\item We develop local and global pipelines for applying magnitude in practice and demonstrate their effectiveness on both simulated TME data~\cite{Bull2023} and tissue microarray images.
\item We compare magnitude-based features with classical descriptors such as cell counts, highlighting their improved ability to capture spatial structure.
\end{itemize}

The remainder of the paper is organised as follows. In Section~2, we describe multispecies spatial data, the ABM model~\cite{Bull2023} that simulates the behaviour of different cell types in the TME, and the spatial proteomics tissue microarray dataset. Section~\ref{sec:magnitude} introduces magnitude and its basic properties. In Section~\ref{sec:results}, we present our feature constructions and experimental results. We conclude with a discussion in Section~\ref{sec:discussion}.

\section{Multispecies Spatial Data}
\label{sec:data}

We refer to point cloud data $X$ whose points carry labels in addition to attributes such as spatial location as \emph{multispecies spatial data}. Specifically, if $X$ has $m$ different species, then $X = \bigcup_{i=1}^{m} X_i$ and each point $x \in X$ is in $\mathbb{R}^n$. Such data arises naturally in applications from spatial biology \cite{carstens2024spatial,jacquemet2020cell, karimi2024method,lewis2021spatial,moffitt2022emerging, park2022spatial, vandereyken2023methods, xiaowei2021method}, but also many other contexts such as geography and ecology, where typically $n=2$.

\subsection{Synthetic data} \label{sec:ABM-data-description}
We use dynamic and spatially-resolved synthetic data generated from an agent-based model (ABM) that simulates the behaviour of different cell types in a TME developed by Bull \& Byrne 2023 \cite{Bull2023} inspired by \emph{in vivo} observations. ABMs are well suited for generating multispecies data because they simulate the emergent behaviour of complex systems through rules governing interactions between individual `agents'. Typical agents are e.g. individual cells \cite{bonabeau2002agent}.

\begin{figure}[h]
  \centering
  \label{fig:a}\includegraphics[width = \textwidth]{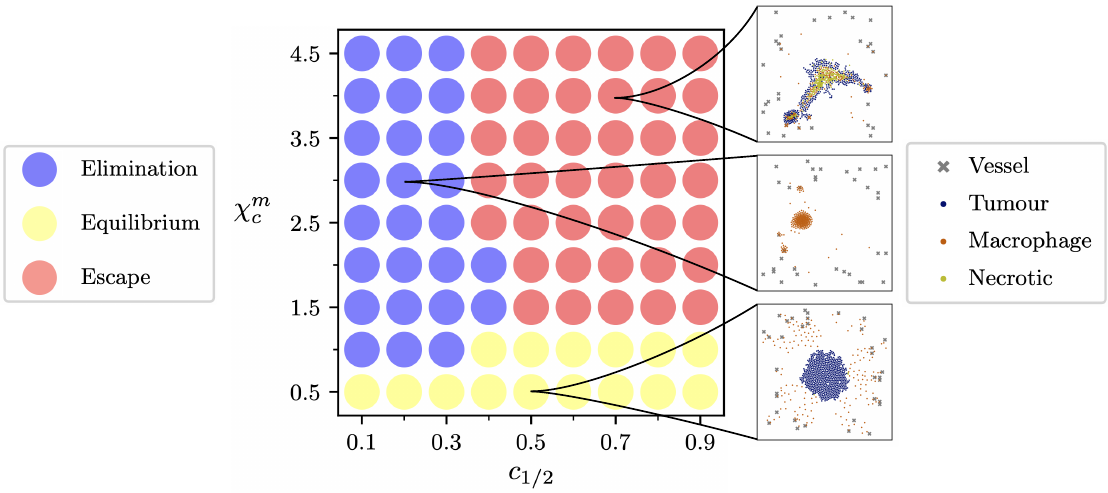}
  \caption{We test our approach on synthetic data generated from an agent based model of the tumour microenvironment (Bull \& Byrne 2023~\cite{Bull2023}). Different combinations of parameter values $c_{1/2}$ (half-maximal macrophage extravasation rate) and $\chi^{m}_c$ (chemotactic sensitivity of macrophages to CSF-1) lead to qualitatively different outcomes: elimination of tumour cells by macrophages (blue parameter regime), equilibrium between tumour cells and macrophages (yellow parameter regime), and escape of tumour cells towards blood vessels (red parameter regime). The pictured classification of qualitative regimes corresponds to the (subjective) classification in the original paper.}
  \label{fig:data-description}
\end{figure}

Each run of the ABM produces a point cloud $X$ with $m = 5$ different species: tumour cells ($X_T$), stromal cells ($X_S$), necrotic cells ($X_N$), macrophages ($X_M$), and blood vessels ($X_V$). With the exception of the blood vessels, whose positions are chosen at random and fixed in the beginning of the simulation, the positions of the different cells change throughout simulated time and are determined by local interactions with other cells and diffusable chemicals (oxygen, CSF-1, TGF\nobreakdash-$\beta$, CXCL12, and EGF). Following Stolz \textit{et al.}\ 2024 \cite{Stolz2024}, we use simulation outputs varying two parameters associated with macrophage behaviour: $\chi^m_c$, the chemotactic sensitivity of macrophages to spatial gradients of one of the chemical species (CSF-1), and $c_{1/2}$, a parameter regulating the rate at which macrophage extravasate from the blood vessels.
Depending on the specific combination of model parameters that are chosen, the simulation outputs cluster naturally into one of three qualitative behaviours that exhibit distinct spatial patterns: elimination of the tumour by the macrophages, equilibrium between the tumour cells and macrophages, and escape of the tumour cells towards blood vessels (see Fig. \ref{fig:data-description}). We use simulation outputs from a predetermined fixed time point (500 hours) where the distinct qualitative behaviours are apparent. We further consider 9 different values for each parameter. For each of the 81 possible parameter pairs $(\chi^m_c, c_{1/2})$, we have up to 20 realisations of the ABM in which the positions of the blood vessels are varied. 
For further details see Stolz \emph{et al.}\ 2024 \cite{Stolz2024}.

\subsection{Colorectal Cancer Data}\label{sec:bio-data-description}
\begin{figure}[h]
  \centering
  \label{fig:a}\includegraphics[width = \textwidth]{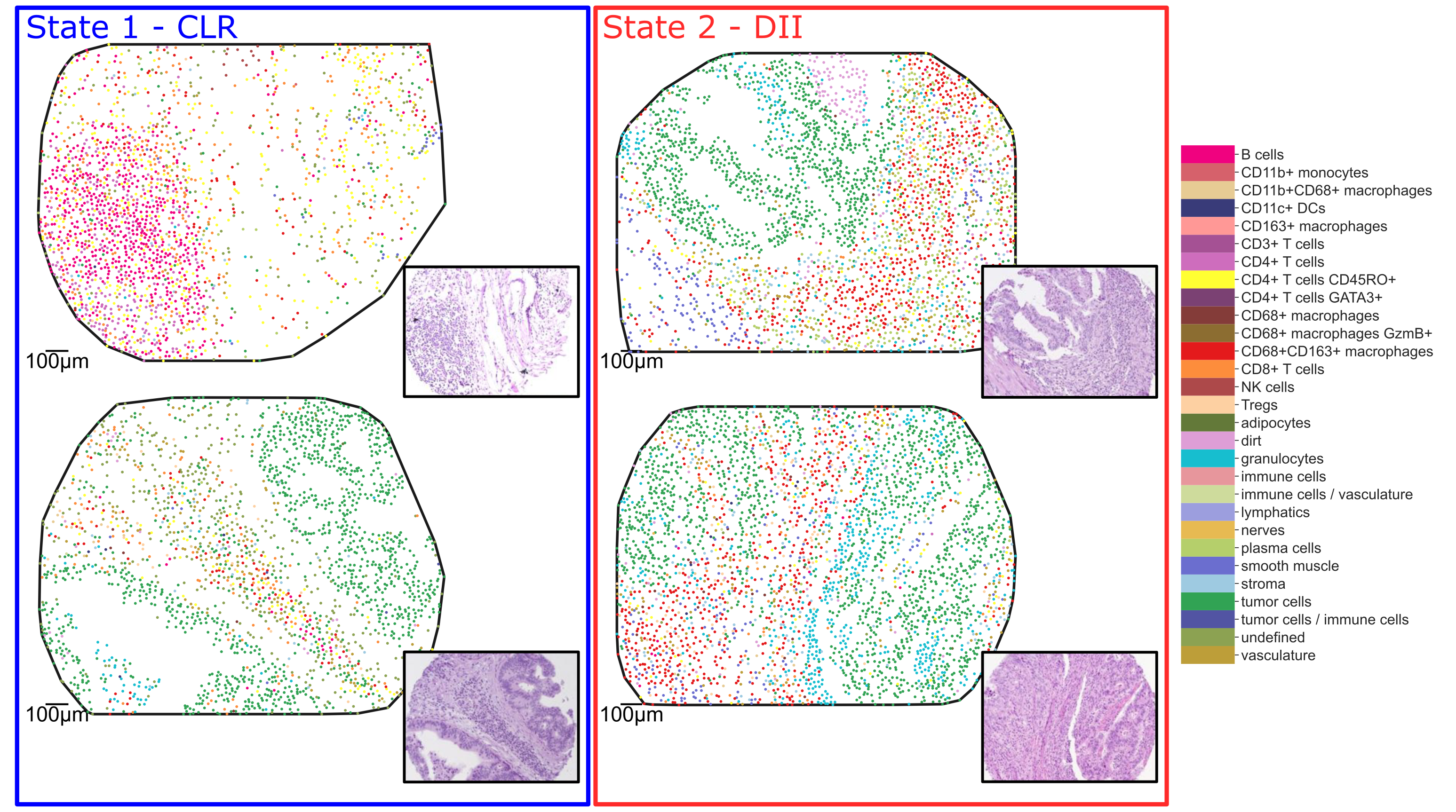}
  \caption{Example point clouds from the colorectal cancer tissue microarray dataset~\cite{schurch_coordinated_2020}, which we use to test our approach. Inset images show the H\&E sections corresponding to these data (reproduced from Sch{\"u}rch \textit{et al.}\ 2020 \cite{schurch_coordinated_2020}). Patients assessed histologically as having a Crohn's like reaction (CLR) were associated with better prognosis than those assessed with diffuse inflammatory infiltration; these assessments reflect opposing ends of a spectrum of immune TMEs.}
  \label{fig:data-description-Nolan}
\end{figure}

We consider 140 tissue microarray (TMA) samples of human colorectal cancer, published by Sch{\"u}rch \textit{et al.}\ 2020~\cite{schurch_coordinated_2020} (see Fig.~\ref{fig:data-description-Nolan}). The samples represent two distinct ends of a spectrum of immune TMEs identified by the authors: `Crohn's like reaction' (CLR) and `diffuse inflammatory infiltration' (DII). CLR is characterised by the presence of tertiary lymphoid structures (clusters of B and T cells that act as local training sites for immune cells), and is associated with good patient prognosis. Conversely, DII is characterised by widespread diffuse (non-structured) infiltration of a wide range of inflammatory immune cell subtypes. The dataset contains a wide range of specific immune cell subtypes, as well as broad labels for other cells present within the TME, with detailed cell classifications and assignments conducted by the original authors.

\section{Mathematical Preliminaries}
\label{sec:magnitude}

In this section we introduce magnitude, a cardinality-like invariant of metric spaces.
It first appeared as the \emph{effective number of species} in ecology \cite{Solow1994} and was later formalised by Leinster in a category theory framework \cite{Leinster2008EulerCharacteristic}.
We restrict attention to finite subsets of Euclidean space, which are most relevant for spatial data.

\subsection{Magnitude and Magnitude Function}\label{sec:magnitude-definitions-properties}

\begin{definition}\label{def:mag}
Let $(P,d)$ be a finite metric space. The \textbf{similarity matrix} of $P$ is
\[
Z_P = \big(e^{-d(p,q)}\big)_{p,q \in P}.
\]
If $Z_P$ is invertible, the \textbf{magnitude} of $P$ is defined by
\[
|P| := \sum_{p,q \in P} (Z_P^{-1})_{p,q}.
\]
\end{definition}
\begin{remark}
More generally, magnitude can be defined for any finite metric space whose similarity matrix admits a weighting~\cite{Leinster2008EulerCharacteristic}, i.e., a vector $w$ satisfying $Z_P w = \mathbf{1}$, even if $Z_P$ is not invertible.
\end{remark}

For finite subsets of $\RR^n$ equipped with the Euclidean metric, the similarity matrix is invertible~\cite[Theorem 2.5.3]{Leinster2008EulerCharacteristic}. In particular, their magnitude is well defined and can be computed as in Definition~\ref{def:mag}.

Magnitude is often interpreted as an \emph{effective number of points}.
Unlike cardinality, it incorporates pairwise distances and thus reflects spatial redundancy: points that lie close together contribute less than well-separated points. The following example, taken from~\cite[Examples 2.1.1, ii]{Leinster2008EulerCharacteristic}, illustrates this interpretation.

    \begin{example}[Two point space]\label{ex.:magnitude-two-point-space}
        Let $P$ be a metric space that consists of just two points $p$ and $q$ with distance $
        \delta$. The similarity matrix \(Z_P\) of \(P\) is
        \[Z_P = \begin{pmatrix}
            1 & e^{-\delta} \\
            e^{-\delta} & 1
        \end{pmatrix},\]
        and the inverse of \(Z_P\) is
        \[Z_P^{-1} = \frac1{1-e^{-2\delta}}\begin{pmatrix}
            1 & -e^{-\delta} \\
            -e^{-\delta} & 1
        \end{pmatrix}.\]
        Thus, the two point space has magnitude
        \[|P| = \frac{2 - 2e^{-\delta}}{1-e^{-2\delta}} = \frac{2}{1+e^{-\delta}}.\]
        Varying the distance $\delta$, we get a continuous, monotonically increasing function \(|P|(\delta)\) on \((0,\infty)\) with values in \((1,2)\) and
        \[\lim_{\delta \to 0} |P|(\delta) = 1, \qquad \lim_{\delta \to \infty} |P|(\delta) = 2,\]
        as demonstrated in Figure~\ref{fig.:two-point-space}. If $p$ and $q$ lie very close to each other, they behave effectively like a single point.
By contrast, if the distance between $p$ and $q$ is large, then they are well separated and thus effectively form two distinct points. It is in this sense that magnitude is interpreted as an effective number of points.
        \begin{figure}[h!]
            \centering
           \includegraphics[scale=0.8]{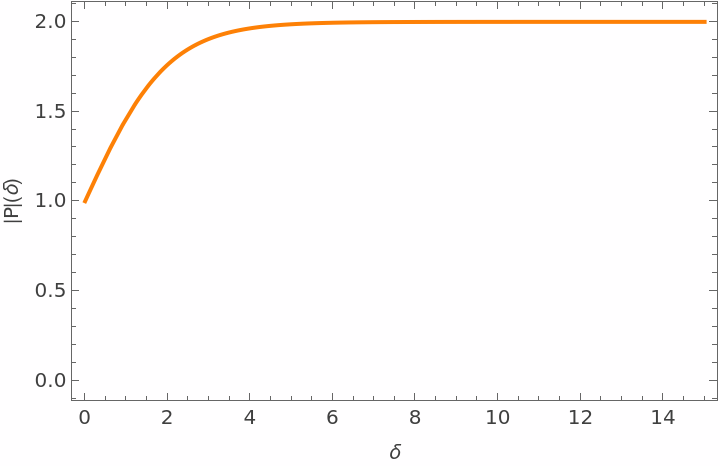}
            \caption{The magnitude \(|P|\) of the two point space as a function of the distance \(\delta\) between the two points.}
            \label{fig.:two-point-space}
        \end{figure}

    \end{example}

Often we are not interested only in the magnitude of $P$, but the magnitude of all its rescaled versions as well. This leads to the notion of the magnitude function~\cite{Leinster2008EulerCharacteristic}.

    \begin{definition}
        Let \(P = (P,d)\) be a finite metric space. For \(t \in (0,\infty)\) denote by \(tP\) the metric space \((P,d_t)\) with \(d_t(p,q) = t\cdot d(p,q)\). Then, the \textbf{magnitude function} of \(P\) is the partially defined function given by
        \[\quad t \mapsto |tP|,\]
        for \(t\in (0,\infty)\) where \(|tP|\) is defined.
    \end{definition}

The magnitude function of subsets of Euclidean space is defined for all \(t\in (0,\infty)\) \cite[Theorem 2.5.3]{Leinster2008EulerCharacteristic}, and it has some very useful properties.

\begin{proposition}\label{prop:properties-of-magnitude}
Let $P \subset \RR^n$ be a finite set equipped with the Euclidean metric.

\begin{enumerate}
    \item \cite[
Corollary 2.4.4]{Leinster2008EulerCharacteristic} If $P \subseteq Q \subset \RR^n$ are finite subspaces of $\RR^n$, then
    \[
    |P| \le |Q|.
    \]
    \item \cite[
Corollary 2.4.5]{Leinster2008EulerCharacteristic} If $P \neq \emptyset$, then $|P| \ge 1$.

    \item \cite[Proposition 2.2.6]{Leinster2008EulerCharacteristic} The magnitude function
    \[
    (0,\infty) \to \RR,
    \qquad t \mapsto |tP|
    \]
    is continuous and satisfies
    \[
    \lim_{t \to \infty} |tP| = \#P,
    \]
    where $\#P$ denotes the cardinality of $P$.
\end{enumerate}
\end{proposition}\

\subsection{Magnitude versus Cardinality}\label{subsec:magnitude_vs_cardinality}

Interpreting magnitude of a finite metric space as an effective number of points naturally raises the question of what the differences between the magnitude of a finite metric space (as an \emph{effective} number of points) and the cardinality of the underlying point set (as an \emph{actual} number of points) are.

While the cardinality of finite sets satisfies the \emph{inclusion-exclusion} principle
\[\#(A \cup B) = \#A + \#B - \#(A \cap B),\]
for all finite point sets $A$ and $B$ (not necessarily disjoint), the same is not true for magnitude.

\begin{figure}[htbp]
  \centering
  \includegraphics[width = \textwidth]{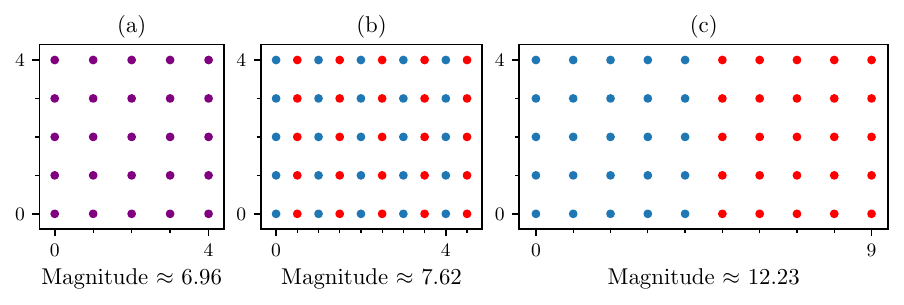}
  \caption{Consider two $5 \times 5$ grids of the scale shown in (a). When the two grids are heavily overlapping, such as in (b), their combined magnitude is not significantly higher than the magnitude of just one grid. When the two grids occupy completely separate spaces, such as in (c), their combined magnitude is almost as large as the sum of the two individual magnitudes. More generally, the difference $|\text{Grid 1}| + |\text{Grid 2}| - |\text{Grid 1} \sqcup \text{Grid 2}|$ (relative to the individual magnitudes) contains information on the spatial relation between the two grids, and in particular on the amount of “shared space” between the two.}
    \label{fig:feature-vec-ex}
 \end{figure}

This failure of inclusion--exclusion is not a drawback but a feature:
the discrepancy
\[
|A| + |B| - |A \sqcup B|
\]
captures the extent to which the two point sets interact spatially as demonstrated in the examples depicted in Figure~\ref{fig:feature-vec-ex}.
This observation motivates the feature constructions introduced in Section~\ref{sec:results}.

\section{Applications}
\label{sec:results}

\begin{figure}[h!]
  \centering
  \label{fig:a}\includegraphics[width = \textwidth]{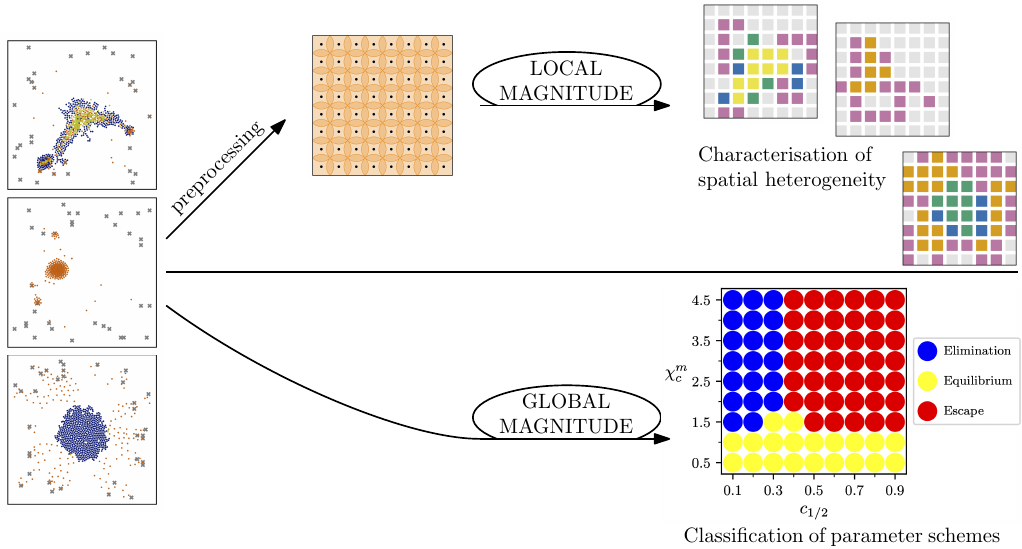}
  \caption{Global and local magnitude pipeline. Input for both pipelines are $20 \times 9 \times 9$ simulated images from different parameter regimes of an agent-based model of the tumour microenvironment (TME). We apply magnitude is globally to generate input for a clustering algorithm that leads to a classification of the parameter regimes into different long term simulation outcomes. After a preprocessing step, we apply magnitude locally to smaller neighbourhoods in order to detect and characterise local architecture in the TME.}
  \label{fig:pipeline}
\end{figure}

We apply magnitude to synthetic multispecies spatial data of the TME from Bull \& Byrne 2023 ~\cite{Bull2023} described in Section~\ref{sec:ABM-data-description}. We do so in two different manners: We gain insight into the association between TME architecture and long-term tumour outcomes by applying magnitude locally; and by applying magnitude globally, we reproduce classification results from Stolz \textit{et al.}\ 2024~\cite{Stolz2024} of different parameter schemes (see Figure~\ref{fig:pipeline}).

In particular, for both applications we have simulation outputs in the form of point sets $X$ in a square domain $[0,50]^2$, together with a partition into different cell types $X_T$, $X_M$, $X_N$, $X_V$ and $X_S$. In our synthetic data stromal cells can be imagined to fill up the entire space not occupied by any of the other cell types and therefore we ignore them in our analysis. Note that this does not necessarily hold in other biological contexts. Further note that as a finite subspace of $\mathbb{R}^2$, these simulation outputs (or subsets thereof) always have magnitude and satisfy the properties outlined in Section~\ref{sec:magnitude} (see, for example, Proposition~\ref{prop:properties-of-magnitude}). For both the local and global application, all magnitudes are taken with a (re-)scaling factor of 0.35. For this particular data and tasks, this scale appears to lie within a small low-scale interval (±0.2) that consistently performs well. For larger scales, where magnitude resembles the number of distinct points, the advantage of taking magnitudes over cell counts is lost as expected.

While for the local application certain “naive” magnitude features prove to be sufficient, we explore more intricate magnitude features in the global application.

Beyond our main application outlined above, we apply our method exemplary to the colorectal cancer data described in Section~\ref{sec:bio-data-description}, demonstrating its applicability to real-world data.

\subsection{Local magnitude}\label{sec:local-magnitude}
\subsubsection{Method - Synthetic Data}
For the local application, we cover the square domain with local neighbourhoods and extract a magnitude feature vector for each of them, in order to then further create local neighbourhood classes.
In particular, we take $8\times8$ closed disks that are evenly spread and cover the domain, as shown in the pipeline in Figure~\ref{fig:pipeline} as the “preprocessing” step. Specifically, we work with a coarse grid $\big\{\frac{50}{2\cdot8},3\frac{50}{2\cdot8},5\frac{50}{2\cdot8},...,15\frac{50}{2\cdot8}\big\}^2$ as centres with radii equal to $\sqrt{2}\frac{50}{2\cdot8}$.
For each such local neighbourhood we take the magnitude $|X_T|$ of tumour cells, the magnitude $|X_M|$ of macrophages and the magnitude $|X_N|$ of necrotic cells within that neighbourhood, yielding a total of $20\times9\times9\times8\times8$ local magnitude feature vectors $(|X_T|,|X_M|,|X_N|)$. We apply $k$-means clustering for $k = 6$ to the resulting magnitude feature vectors in order to obtain a classification of the local neighbourhoods into different local classes. Here, we use the implementation from scikit-learn~\cite{scikit-learn}. We deem the choice $k=6$ suitable from a subjective evaluation of the resulting classification. For an analysis of different values for $k$ and corresponding results, we refer to Appendix~\ref{app:number-of-clusters}.

\subsubsection{Local Magnitude Detects Spatial Heterogeneity in the Simulated TME}
\begin{figure}[htbp]
  \centering
  \label{fig:a}\includegraphics{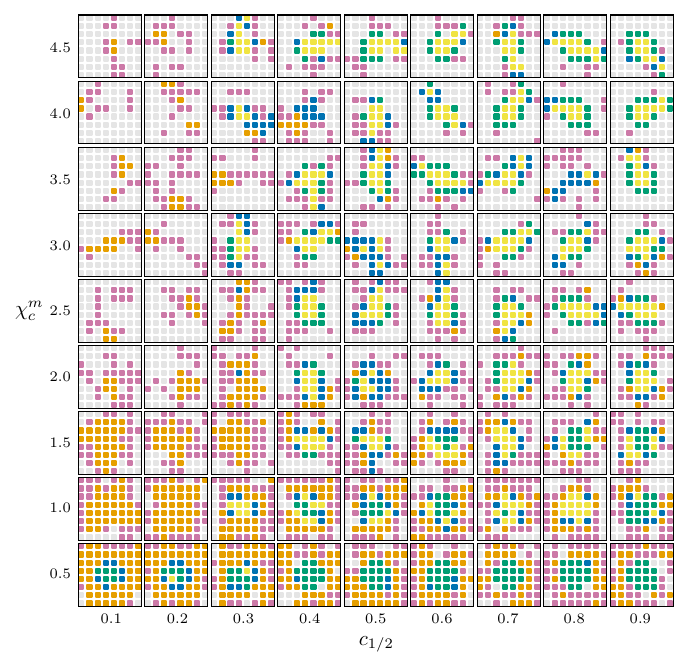}
  \caption{Categorisation of local neighbourhoods into 6 classes according to the local magnitude features. Depicted here are for each parameter scheme the labelled neighbourhoods of one exemplary simulation outcome from the ABM~\cite{Bull2023}.}
  \label{fig:local-across-schemes}
\end{figure}

The results of the experiments described above are $8\times 8$ neighbourhoods per simulation outcome, each assigned to one of 6 classes. Figure~\ref{fig:local-across-schemes} depicts these outcomes for one exemplary and pseudo-randomly chosen simulation outcome per parameter scheme. We use these same $9 \times 9$ representative simulation outcomes (or a $5 \times 5$ subset thereof) throughout this work, and the simulation outcomes themselves are pictured for reference in Figure~\ref{fig:plots-across-parameter-schemes} in the appendix.

\begin{figure}[htbp]
  \centering
  \label{fig:a}\includegraphics{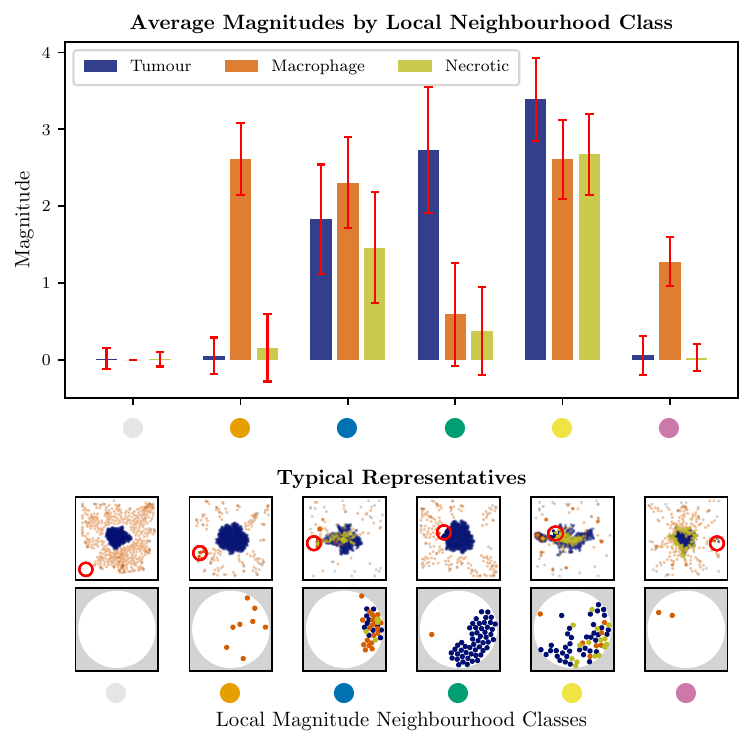}
  \caption{Mean features of each local magnitude neighbourhood class, together with a typical representative for each class. Each local magnitude neighbourhood class is labelled by a coloured disk, corresponding to the same coloured classes in Figure~\ref{fig:local-across-schemes}.}
  \label{fig:signatures}
\end{figure}

Most of the simulation outcomes show a clear radial pattern in their local magnitude class, and subjectively, one can recover the same qualitative classification of parameter schemes into long term qualitative outcomes of the simulations as described by Bull and Byrne in~\cite{Bull2023} and Stolz \textit{et al.}\ 2024~\cite{Stolz2024} (see Figures~\ref{fig:results-bull},~\ref{fig:results-stolz}) – results that will be reproduced and discussed in the global application in Section~\ref{sec:global-magnitude}. Note that different seeds for the pseudo-random choice of displayed results give the same qualitative observations. Further, we consider the mean magnitude features, i.e., the cluster centres, of each of the six neighbourhood classes. We recover six distinct local neighbourhood classes that we qualitatively describe as the “empty neighbourhood” (\textcolor[HTML]{E5E5E5}{\Large $\bullet$}), “high macrophage neighbourhood” (\textcolor[HTML]{E69F00}{\Large $\bullet$}), the “low mixed neighbourhood” (\textcolor[HTML]{0072B2}{\Large $\bullet$}), the “tumour neighbourhood” (\textcolor[HTML]{009E73}{\Large $\bullet$}), the “high mixed neighbourhood” (\textcolor[HTML]{F0E442}{\Large $\bullet$}), and the “low macrophage neighbourhood” (\textcolor[HTML]{CC79A7}{\Large $\bullet$}), with terms “high” and “low” referring to the respective mean magnitude. The mean features of each neighbourhood class, as well as a typical representative, are depicted in Figure~\ref{fig:signatures}. The representatives were chosen by selecting the one simulation outcome with magnitude feature vector closest to the mean features (with respect to the Euclidean metric).

\begin{figure}[htbp]
  \centering
  \label{fig:a}\includegraphics{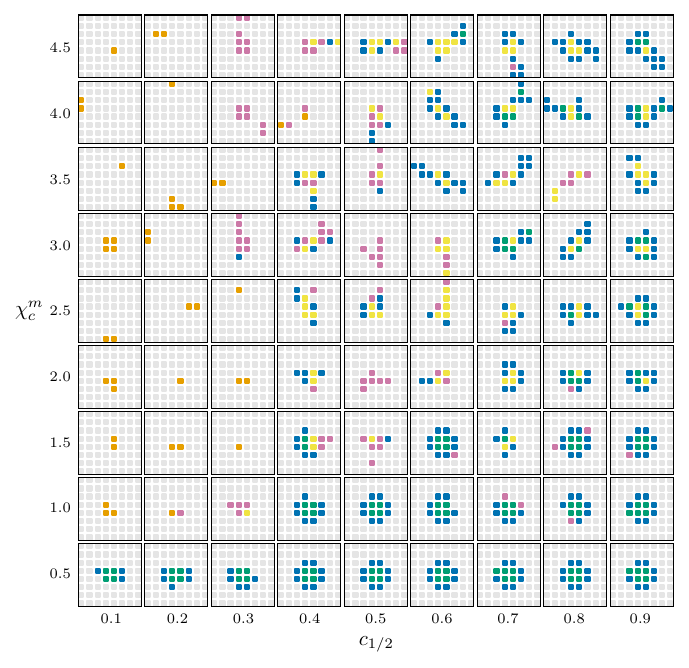}
  \caption{Categorisation of local neighbourhoods into 6 classes according to the local cell counts. Beyond the replacement of magnitude with cell counts, the exact method from Section~\ref{sec:results} is used. Represented here are the same exemplary simulation outcomes from the ABM~\cite{Bull2023} as in Figure~\ref{fig:local-across-schemes}.}
  \label{fig:cellcounts}
\end{figure}

As discussed in Section~\ref{sec:magnitude}, magnitude can often be interpreted as an effective number of points, or, in the context here, as a “cell count plus”, which takes into account the spatial distribution of cells. We compare the use of naive magnitude features with cell counts by replacing magnitudes with cell counts in the described method. The results in Figure~\ref{fig:cellcounts} confirm that cell counts alone detect significantly less spatial heterogeneity in the TME. Cell counts detect a limited radial pattern in the densest region of the TME but appear to be blind towards local architecture in more sparse regions. In particular, the naive approach with cell counts fails to detect different spatial distribution – and thus behaviours – of macrophages. We emphasise that the behaviour of macrophages is particularly relevant for the model and correlates heavily with long term qualitative simulation outcomes as described in~\cite{Bull2023}.
Further analysis shows that only when combined with average pairwise distances per cell type and normalisation prior to clustering, cell count feature vectors yield \textit{almost} comparable results. However, those results still contain less of the peripheral macrophage structure compared to the magnitude results in Figure~\ref{fig:local-across-schemes}. Supporting figures and more details can be found in Appendix~\ref{app:cellcounts}.

Comparing magnitude to cell counts can provide insight into how magnitude behaves differently and, consequently, into when each type of feature is most appropriate. Most notably, from Figures~\ref{fig:local-across-schemes} and~\ref{fig:cellcounts} we observe that magnitude enhances architecture. Boundary regions of cell clusters may have significantly fewer cells than the dense centers of such clusters, but the cells may be spread out sufficiently for magnitude features to still detect them as part of the cluster. Consequently, magnitude can pick up on more subtle radial patterns with a coarser neighbourhood structure. Consider, for example, the subtle ring of necrotic cells that appears frequently for simulation outcomes at $\chi_c^m \in \{1,1.5\}$ (see Figure~\ref{fig:plots-across-parameter-schemes}). Magnitude-based features in Figure~\ref{fig:local-across-schemes} detect this, as can be observed most notably in the subplots at parameter combinations $(c_{1/2},\chi_c^m) \in \{(0.4,1),(0.5,1),(0.6,1),(0.7,1)\}$. There, yellow, necrotic-dense neighbourhoods cluster around a small, green, tumour-dense core rather than fully constituting the center themselves.

\subsubsection{Method - TMA Data}
Following our approach for the synthetic data, we cover each TMA sample spot with a lattice of local neighbourhoods. Here, we work with neighbourhoods with radius $r=200\mum{}$, centred on a coarse grid with coordinates $\big\{min(r)+100, min(r)+300, ...,max(r)-100 \big\}^2$ (where $min(r)$ and $max(r)$ represent the extreme $x$ or $y$ coordinates of the bounding box). For each local neighbourhood, we take the magnitude $|X_i|$ of each of the 27 possible cell types (excluding `dirt' and `undefined'), where $i \in \{0,...,26\}$ indexes the cell type. As before, we apply $k$-means clustering with $k=6$ to the resulting magnitude feature vectors (with $k$ determined by subjective evaluation).

\subsubsection{Local Magnitude Clusters Differentiate CLR from DII}
\begin{figure}[htbp]
  \centering
  \includegraphics[width=.7\textwidth]{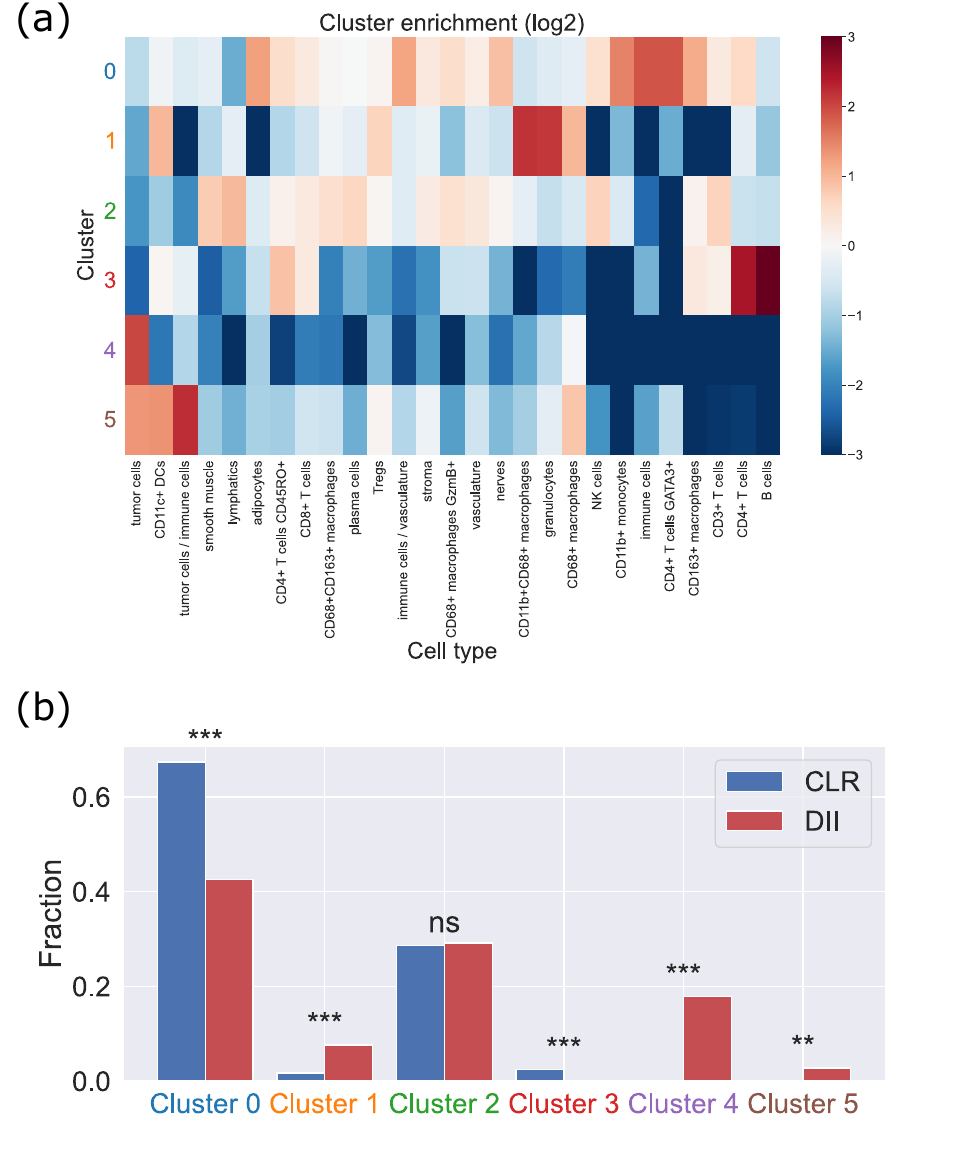}
  \caption{(a) Cluster enrichment scores for cell types within local neighbourhoods defined by magnitude. Cluster 3 corresponds to the `follicle' cellular neighbourhood defined in~\cite{schurch_coordinated_2020}, which was found to be significantly associated with CLR over DII.\\
  (b) In our analysis, Cluster 3 is still significantly associated with CLR over DII. We find additional clusters which are over-represented in CLR (Cluster 0) or DII (Clusters 1, 4, 5). $p$-values were calculated using Fisher's exact test with Benjamini-Hochberg correction ($*: p<0.05$, $**: p<0.01$, $***: p<0.001$).}
  \label{fig:nolan_local}
\end{figure}

After clustering neighbourhoods based on local magnitude, we characterise them according to their enrichment by the magnitude of their constituent cell types (Fig.~\ref{fig:nolan_local}(a)). Critically, Cluster 3 in our analysis has a similar profile to the `follicle' cellular neighbourhood identified by Sch{\"u}rch \textit{et al.}\ 2020~\cite{schurch_coordinated_2020} (enriched for B cells and CD4+ T cells, with lower enrichment for CD163+ macrophages and CD3+ T cells). Our analysis also finds this neighbourhood over-represented in CLR patients compared to DII patients (Fig.~\ref{fig:nolan_local}(b)), reinforcing the observations made in the original study. However, defining the neighbourhoods using magnitude reveals additional differences between the conditions which the original, cell-count based, cellular neighbourhoods were unable to identify. In particular, Cluster 4 and 5 are substantially overrepresented in patients with DII (poor prognosis); these neighbourhoods have high magnitudes of tumour cells and tumour associated immune cells. Cluster 1 is also associated with worse prognosis, and is characterised by high magnitudes of infiltrating immune cells including CD68+ macrophages, CD11b+CD68+macrophages, and granulocytes.

\subsection{Global Magnitude}
\label{sec:global-magnitude}

Applying magnitude to the synthetic data \textit{globally}, i.e., to the entire simulation domains, we reproduce results from the work of Bull and Byrne~\cite{Bull2023} and Stolz \textit{et al.}\ 2024~\cite{Stolz2024} (see Appendix~\ref{app:global-previous}). We use the simulation outcomes~\cite{Bull2023} from the $9 \times 9$ parameter schemes to classify the parameter schemes into different qualitative behaviours and compare these to known ground truth~\cite{Bull2023,Stolz2024}. We explore and compare different types of magnitude features.

\subsubsection{Method - Synthetic Data}
We compute magnitude-based feature vectors for all simulation outcomes, and then proceed with the classification as in Stolz \textit{et al.}\ 2024~\cite{Stolz2024}:
The resulting feature vectors are clustered using $k$-means clustering for $k=3$. This yields a partition of all simulations into the three qualitative behaviours Elimination, Equilibrium and Escape (see Fig.~\ref{fig:data-description}). For each parameter combination, we depict the most dominant cluster across 20 simulations and calculate the \emph{purity} of each parameter scheme, i.e., the number of outcomes assigned to the dominant behaviour relative to the total number of outcomes within the given parameter scheme.

Here, we extend the simple magnitude feature vectors $(|X_T|, |X_M|, |X_N|)$ that proved sufficient in our local magnitude application by adding \textit{inclusion-exclusion type differences} of magnitudes (as already hinted at in Section~\ref{subsec:magnitude_vs_cardinality}).
Given a multispecies point set \(X = X_1 \sqcup \cdots \sqcup X_k\), we calculate differences
\[|X_{i_1}| + \dots + |X_{i_m}| - |X_{i_1}\sqcup \dots \sqcup X_{i_m}|\]
for non-empty subsets \(\{i_1,...,i_m\}\) of \([k] := \{1,...,k\}\).
All such differences together yield a feature vector
\begin{equation}\label{eq:incl-excl}
    \left(\sum_{i \in I}|X_i|  - \left|\bigsqcup_{i\in I}X_i\right| \quad \text{for } I \in \big(2^{[k]}\setminus\{\emptyset\}\big) \right).
\end{equation}

For our data, we concatenate the simple feature vector $(|X_T|,|X_M|,|X_N|,|X_V|)$ with the feature vector described above with species $X_T$, $X_M$, $X_N$ and $X_V$:

\begin{equation}\label{eq:incl-excl-explicit}
    \begin{aligned}
    &(|X_T|,|X_M|,|X_N|,|X_V|,\\
     &\quad |X_T|+|X_M| - |X_T \sqcup X_M|,\ \ldots,\\
    &\quad\quad|X_T|+|X_M|+|X_N|+|X_V|-|X_T\sqcup X_M \sqcup X_N \sqcup X_V|)
    \end{aligned}
\end{equation}

In contrast to the local magnitude application, we include blood vessels, as the spatial relation of the blood vessels with other cell types contains information about the state of the TME. Just like in the local magnitude application, these feature vectors include a choice of scale at which the magnitude is taken. Note that one could also consider inclusion-exclusion type differences of the form
\[|X_{i_1}\sqcup...\sqcup X_{i_l}| + |X_{j_1}\sqcup...\sqcup X_{j_m}| - |X_{i_1}\sqcup...\sqcup X_{i_l} \sqcup X_{j_1}\sqcup...\sqcup X_{j_m}|\]
for some nonempty disjoint subsets \(\{i_1,...,i_l\}, \{j_1,...,j_m\}\) of \([k] := \{1,...,k\}\), though taking all such differences would quickly lead to feature vectors that are no longer linearly independent.

\subsubsection{Global Magnitude Confirms Existing Classification of Synthetic Data}
\begin{figure}[htbp]
  \centering
  \label{fig:a}\includegraphics[]{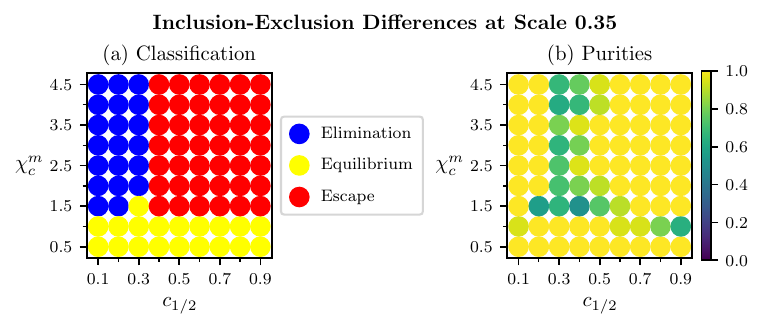}
  \caption{(a) Classification of \(9 \times 9\) pairs of parameter values \(\chi_c^m\) (chemotactic sensitivity of macrophages to CSF-1) and \(c_{1/2}\) (half-maximal macrophage extravasation rate) resulting from clustering magnitude feature vectors that include inclusion-exclusion type differences (as defined in Equation \eqref{eq:incl-excl},~\eqref{eq:incl-excl-explicit}). (b) Purity scores of the classification. That is, for each parameter pair the number of outcomes assigned to the dominant behaviour relative to the total number of outcomes within the given parameter scheme.}
  \label{fig:results-global-magnitude}
\end{figure}

We show the classification of parameter schemes that results from concatenated magnitude feature vectors, with magnitudes taken at scale 0.35, in Figure~\ref{fig:results-global-magnitude}.

\begin{figure}[h!]
  \centering
  \label{fig:a}\includegraphics[]{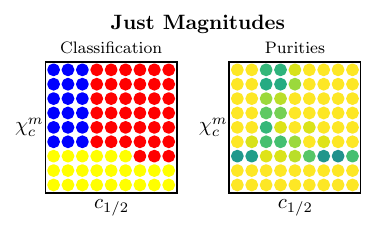}
  \hspace{.5cm}
  \label{fig:b}\includegraphics[]{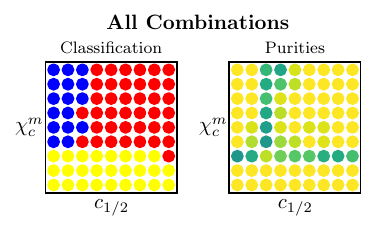}
  \caption{Classification and purity scores of \(9 \times 9\) pairs of parameter values \(\chi_c^m\) (chemotactic sensitivity of macrophages to CSF-1) and \(c_{1/2}\) (half-maximal macrophage extravasation rate) resulting from clustering magnitude feature vectors that consist only of the “naive” features $|X_S|$ for $S \in \{T,M,N,V\}$ (left), or that consist of the magnitudes of all possible combinations of cell types (as defined in~\eqref{eq:all-comb}), but no inclusion-exclusion type differences.}
  \label{fig:results-global-magnitude-comparison}
\end{figure}

Verification and comparison of these results with results of previous methods can be found in Appendix~\ref{app:global-previous}. Most notably, our results agree both with the subjective classification shown in Figure~\ref{fig:results-bull} and with the TDA-based classification in Figure~\ref{fig:results-stolz}. Inclusion-exclusion type feature vectors appear to produce a better classification than simple magnitude feature vectors $(|X_T|,|X_M|,|X_N|)$, whose corresponding result is shown in Figure~\ref{fig:results-global-magnitude-comparison} (left). This holds both when compared against previous results and when assessed through direct visual inspection of the actual simulation outcomes.
For further comparison, Figure~\ref{fig:results-global-magnitude-comparison} (right) shows the classification that results from taking feature vectors consisting of the magnitudes of all possible combinations of species, i.e. $$(|X_T|,...,|X_V|,|X_T\sqcup X_M|,...,|X_T \sqcup X_M \sqcup X_N \sqcup X_V|),$$ or formally:
\begin{equation}\label{eq:all-comb}
    \left(\left|\bigsqcup_{Y\in \mathcal{Y}} Y\right| \quad \text{for } \emptyset \neq \mathcal{Y} \subseteq \{X_T,X_M,X_N,X_V\}\right).
\end{equation}

The authors of~\cite{Stolz2024} also show that these results cannot be reproduced using simple descriptor vectors consisting of cell counts and average distances of each cell type to the nearest blood vessel (for a comparison, see Figure 1 in SI of~\cite{Stolz2024}). Hence, this example indicates that the non-additivity with respect to disjoint unions can be helpful and giving magnitude an advantage over cell counts.

\subsubsection{Global Magnitude for TMA Data Reinforces Analysis from Local Magnitude}
\begin{figure}[h!]
  \centering
  \includegraphics[width=\linewidth]{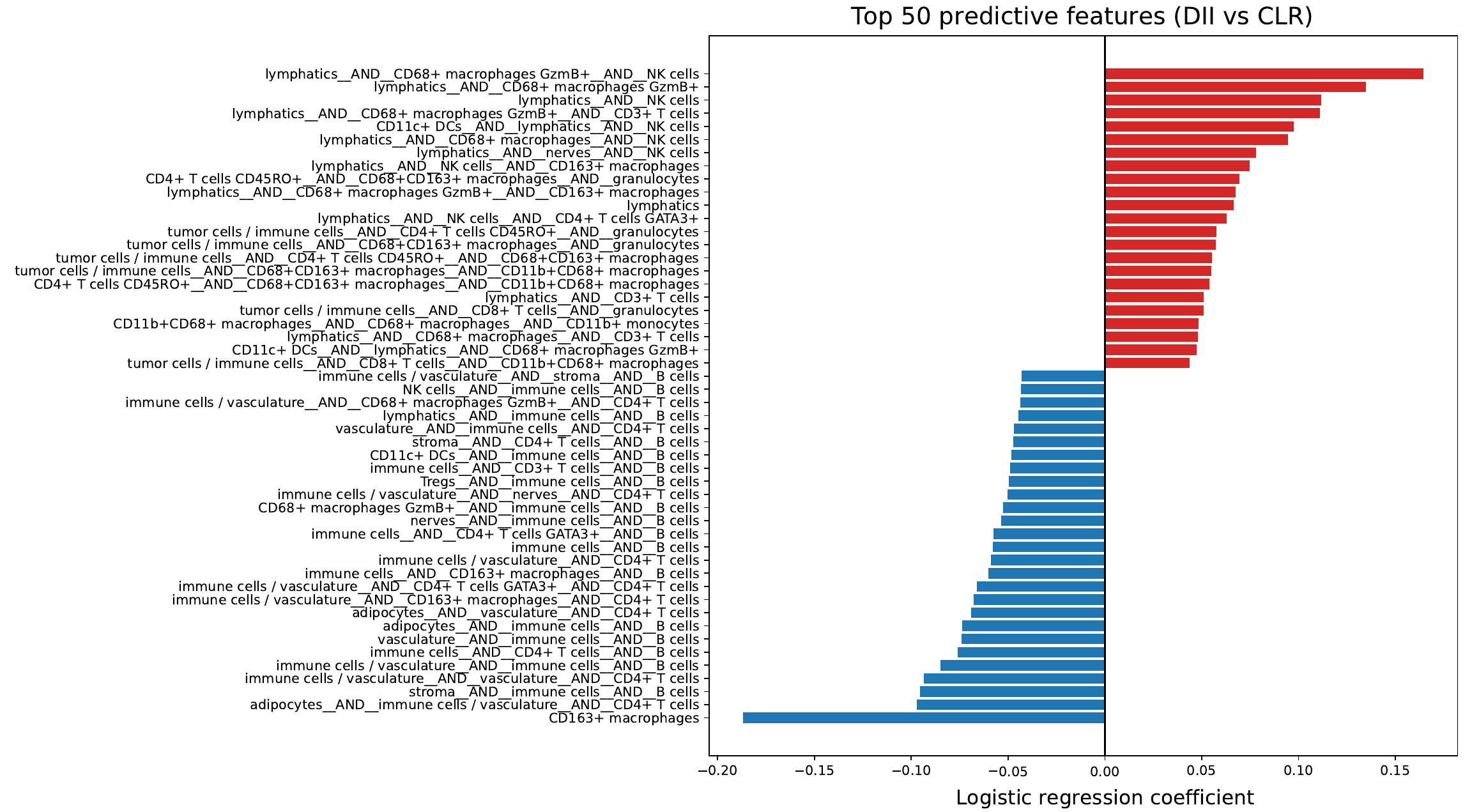}
  \caption{Most important 50 features contributing to model fit in the 5-fold cross-validation, according to logistic regression coefficient. Red features are more strongly associated with DII (poor prognosis) than CLR, while blue features are more strongly associated with good prognosis. The phrase \_\_AND\_\_ in labels indicates an inclusion.}
  \label{fig:results-global-magnitude-comparison-real}
\end{figure}

As with the synthetic data, we also calculate global magnitudes for the TMA dataset (i.e., magnitude across the entire TMA spot, rather than within a local neighbourhood). For each TMA spot, we calculate the global magnitude $|X_i|$ for all cell types $i$, and the combinations $|X_i \sqcup X_j|$ and $|X_i \sqcup X_j \sqcup X_k|$ for all triplets of distinct cell types $i$, $j$ and $k$. Concatenating these generates a 3,303 dimensional feature vector, which we use to identify the magnitudes most strongly discriminatory between CLR and DII. We use scikit-learn~\cite{scikit-learn} to apply a standard scaler and perform logistic regression (\verb|max_iter=1500|), followed by 5-fold cross-validation (AUC: $0.771 \pm 0.116$). The most important 50 features, as identified by their logistic regression coefficient, are shown in Fig.~\ref{fig:results-global-magnitude-comparison-real}.

Reinforcing the results of both our local magnitude study and those reported in \cite{schurch_coordinated_2020}, a high magnitude of B cells, in combination with a broad spectrum of immune cell subtypes, is most associated with CLR and hence good patient prognosis. This suggests that a TME containing B~cell follicles or tertiary lymphoid structures is indeed highly associated with good patient prognosis. The single most important feature is the magnitude of CD163+ macrophages in the TME, which we observed to be enriched in our local analysis in both Cluster 0 and Cluster 3 (see Fig.~\ref{fig:nolan_local}), two clusters associated with CLR. The cell types most associated with DII and poor prognosis are associations between NK cells, lymphatics, and CD68+ GzmB+ macrophages. Notably, when these cell types are present the inclusion of CD163+ macrophages is associated with DII, suggesting that macrophage phenotypic plasticity (influenced by the TME) may play a key role in distinguishing CLR from DII.

\section{Discussion}\label{sec:discussion}

In this manuscript, we used magnitude feature vectors to extract spatial information from multispecies spatial data and demonstrated their potential through applications to synthetic TME data.

Magnitude, in the special case of finite Euclidean subspaces, has particularly favorable properties (see Section~\ref{sec:magnitude-definitions-properties}), supporting an interpretation as an effective number of points. Our results suggest that magnitude holds considerable promise for the analysis of spatial multispecies data. In particular, magnitude captures spatial information in the TME both locally and globally, and shows strong alignment with qualitative long-term simulation outcomes. Local magnitude provides insight into spatial architecture, while global magnitude reproduces known classifications of parameter regimes governing long-term behaviour.

We compare magnitude-based features with simple cell counts, which are a standard summary statistic for multispecies data and coincide with magnitude in the asymptotic large-scale regime. We show that at sufficiently small scales, the additional spatial information captured by magnitude leads to substantially improved performance in data analysis. This improvement is not easily replicated by augmenting cell counts with pairwise distance information.

Quantities used to assess spatial heterogeneity are often scale-dependent \cite{spatialYuan2016}. In contrast, magnitude provides a flexible multiscale framework, as it can be evaluated and tuned across a range of spatial scales. Moreover, magnitude is highly versatile in multispecies data analysis: while clustering based on magnitude values for each cell type is already effective, it also supports more expressive features beyond those derived from cell counts. For example, multiscale information can be incorporated, and the non-additivity of magnitude under disjoint unions can be exploited through inclusion--exclusion type expressions, such as $|X_1| + |X_2| - |X_1 \sqcup X_2|$, which capture interactions between species and provide additional insight into spatial relationships.

Here, we used synthetic data due to its controlled setting and the availability of ground truth. We further validated our work by application to a real dataset of 140 tissue microarray samples from human colorectal cancer, identifying the key relationships highlighted by previous cell-count based analyses while finding additional cell types of interest which are not suggested by cell counts alone. While our approach shows great promise, future work will require validation on a broader range of real datasets. In real data, spatial proximity of cells often serves as a proxy for potential cell–cell interactions. Although our method enables a finer characterisation of cell clusters than approaches based solely on cell counts, it does not directly infer cell–cell communication, which requires additional modelling. Spatial omics data, including imaging-based colocalisation and spatial transcriptomics, can help characterise interactions within niches \cite{birk2025quantitative, bull_muspan_2025, carstens2024spatial, moffitt2022emerging, park2022spatial, vandereyken2023methods, schurch_coordinated_2020}. It would therefore be interesting to combine our features with methods that explicitly model such interactions by integrating multiple modalities, such as NicheCompass \cite{birk2025quantitative}. Further work is also needed to link observed spatial patterns to functional and clinical outcomes \cite{bareham2024defining, carstens2024spatial}.

Existing approaches to multispecies data analysis include summary statistics of counts and locations \cite{goltsev2018deep}, nearest-neighbour methods \cite{bull_muspan_2025, palla2022squidpy, parra2021methods}, spatial statistics \cite{Bull2023, bull_muspan_2025, bull2024extended}, point process models \cite{carstens2017spatial}, neighbourhood-based summaries \cite{bull_muspan_2025, palla2022squidpy}, TDA~\cite{bae2025stopover, di_Montesano_2025, fu2025persistent, natarajan2026topology,park2025spatial, Stolz2024,vipond2021multiparameter, yang2025topological}, and machine learning approaches \cite{studer2023tumor} (see, e.g., \cite{carstens2024spatial} for a detailed list of common approaches). Magnitude complements these methods by providing low-dimensional, interpretable features with a clear geometric meaning, while remaining flexible and extensible.

\section{Code Availability}
The code is available at \url{https://github.com/j-sollberger/tumor-magnitude}.

\bibliographystyle{plain}
\bibliography{references}
\clearpage

\appendix
\section{Further Local Magnitude Variations and Comparisons}
\begin{figure}[!h]
  \centering
  \label{fig:a}\includegraphics{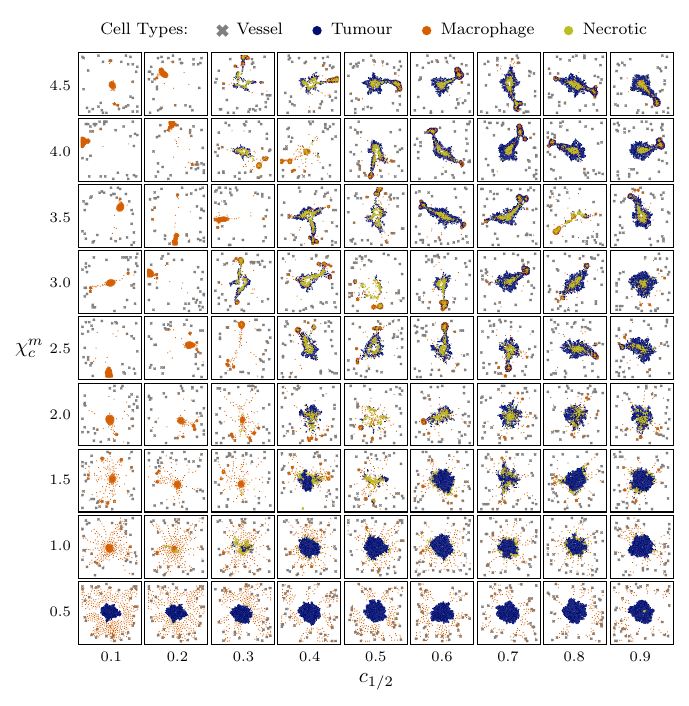}
  \caption{Simulations across parameter schemes for the pseudo-randomly chosen parameter schemes that are used to display results throughout this paper.}
  \label{fig:plots-across-parameter-schemes}
\end{figure}

Further categorisations of local neighbourhoods are given using more intricate non-topological feature vectors using cell counts and pairwise distances, as well as different numbers of clusters for both magnitude and cell count feature vectors. For all local magnitude results that are displayed across parameter schemes (Figures~\ref{fig:local-across-schemes},~\ref{fig:cellcounts} in the main text and Figures~\ref{fig:making-celcounts-work},~\ref{fig:local-magnitudes-5signatures},~\ref{fig:local-magnitudes-7signatures},~\ref{fig:local-magnitudes-9signatures} in the appendix), the same pseudo-randomly selected simulation outcomes are used as representatives for the parameter schemes. In order to better compare these results subjectively, plots of the respective simulation outcomes are given in Figure~\ref{fig:plots-across-parameter-schemes}.

\subsection{More intricate method with cell counts and pairwise distances}
\label{app:cellcounts}
In the local magnitude application, replacing the simple magnitude feature vectors $(|X_T|,|X_M|,|X_N|)$ with simple cell count feature vectors $(\#X_T,\#X_M,\#X_N)$ did not yield comparable results for the method and purpose described in Section~\ref{sec:local-magnitude}.

\begin{figure}[h!]
  \centering
  \label{fig:a}\includegraphics{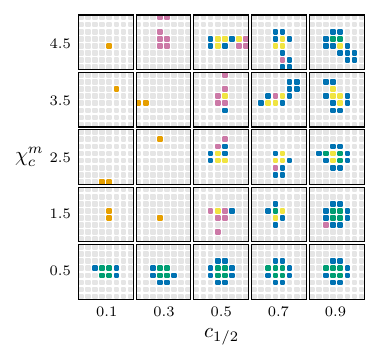}
  \label{fig:b}\includegraphics{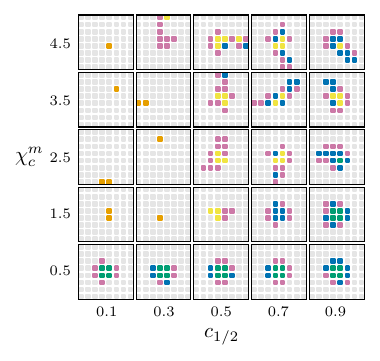}
  \label{fig:c}\includegraphics{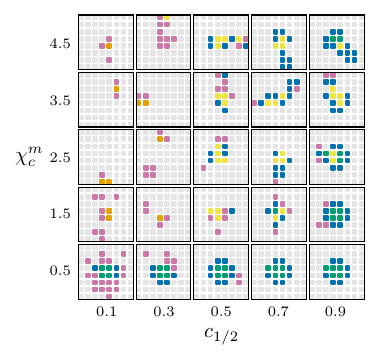}
  \label{fig:d}\includegraphics{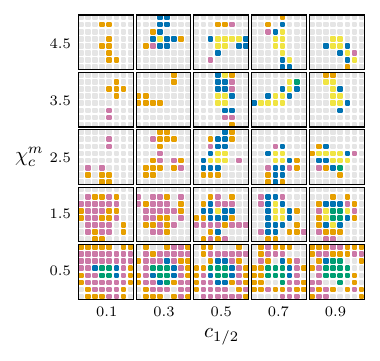}
  \caption{Categorisation of local neighbourhoods into 6 classes according to local cell counts with or without average pairwise distance per cell type and with or without normalisation prior to clustering. Top left: Categorisation according to cell counts without normalisation. Top right: Categorisation according to cell counts with normalisation. Bottom left: Categorisation according to cell counts and distances without normalisation. Bottom right: Categorisation according to cell counts and distances with normalisation.}
  \label{fig:making-celcounts-work}
\end{figure}

We optionally extend the simple cell count feature vectors with the average pairwise distance between each two cell types, and we optionally normalise the cell count feature vectors prior to the $k$-means clustering. For the fixed number of 6 classes, Figure~\ref{fig:making-celcounts-work} shows the categorisations of local neighbourhoods for simple cell count feature vectors (as already done for 5 classes in the main text), for normalised cell count feature vectors, for feature vectors consisting of cell counts as well as average pairwise distance per cell type, and for normalised versions of these same cell counts plus average distances feature vectors. These serve to support the claims in our comparison to cell counts in Section~\ref{sec:results}.

\subsection{Selection of the Number of Clusters}
\label{app:number-of-clusters}
\begin{figure}[htbp]
  \centering
  \label{fig:a}
  \includegraphics{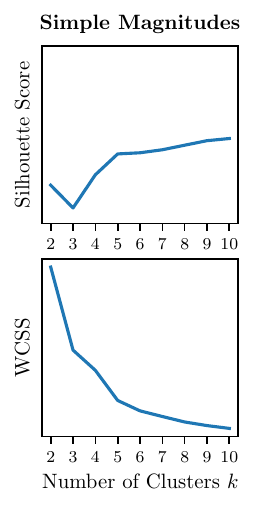}
  \hspace{.5cm}
  \includegraphics{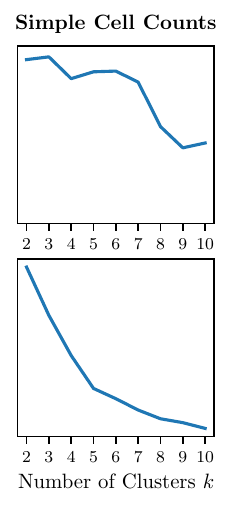}
  \hspace{.5cm}
  \includegraphics{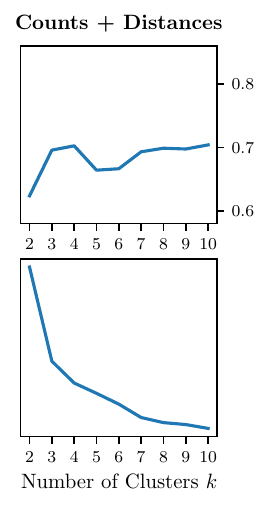}
  \caption{Cluster quality analysis of the $k$-means clustering of local TME neighbourhoods for simple magnitude feature vectors, simple cell count feature vectors, and more intricate normalised feature vectors containing cell counts and average pairwise distances between two cell types.}
  \label{fig:clustering-scores}
\end{figure}

One possibility for variation in our method from Section~\ref{sec:local-magnitude} is the choice of $k$ in the $k$-means clustering.
For the subjective evaluation of results presented in this work, we aim for a cluster number between 4 and 8 in order to ensure a balance between sufficient detail and low complexity suitable for visualisation. Here, we assess the cluster quality for the presented method with simple magnitude feature vectors $(|X_T|,|X_M|,|X_V|)$, for the direct comparison with simple cell count feature vectors $(\#X_T,\#X_M,\#X_V)$, and for the further exploration of standard methods in Appendix~\ref{app:cellcounts} with feature vectors which consist of cell counts as well as average pairwise distances and are normalised. The silhouette score~\cite{ROUSSEEUW_1987} captures a combination of separation between and cohesion within clusters, and may infer whether with increasing values of $k$ more meaningful structure can be discovered. The elbow method with respect to within cluster sums of squares (WCSS) suggests a value for $k$ beyond which overall compactness of clusters does not significantly improve anymore, and thereby also indicates a minimal value for $k$ necessary to avoid missing meaningful structure~\cite{scikit-learn}. For both, the implementation by Scikit-learn~\cite{scikit-learn} was used, and the results are depicted in Figure~\ref{fig:clustering-scores}.
The relatively high silhouette scores between $0.6$ and $0.85$ as well as at least one clear elbow appearing in each WCSS curve indicate that for all three represented choices of features there is some cluster structure present in the data that can be captured using $k$-means.

\begin{figure}[!h]
  \centering  \label{fig:a}\includegraphics{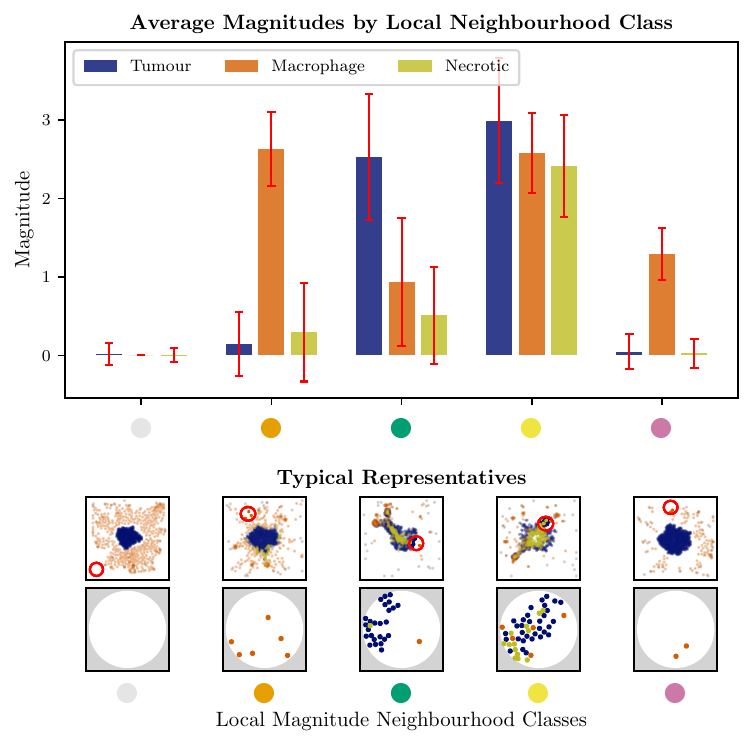}
  \caption{Average features in the 5 classes after categorisation of local neighbourhoods according to the local magnitudes.}
  \label{fig:averages-magnitudes-5signatures}
\end{figure}

\begin{figure}[htbp]
  \centering
  \label{fig:a}\includegraphics{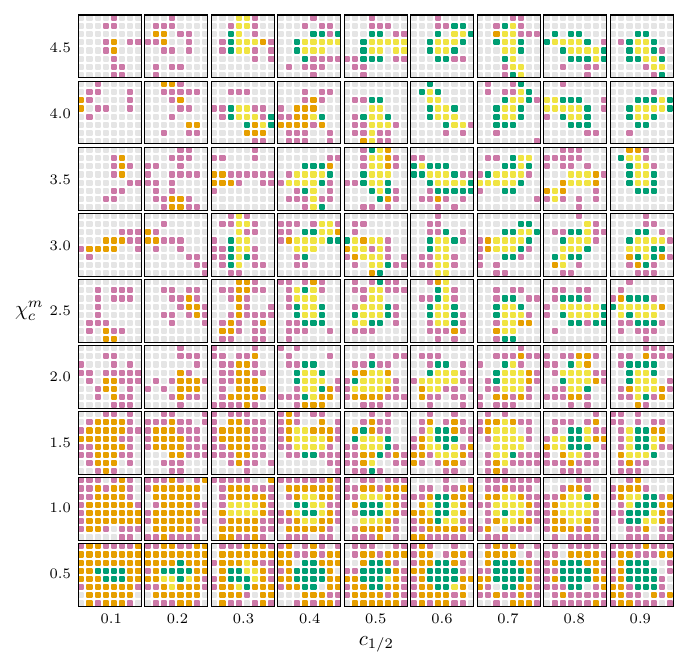}
  \caption{Categorisation of local neighbourhoods into 5 classes according to the local magnitudes of tumour cells, macrophages and necrotic cells.}
  \label{fig:local-magnitudes-5signatures}
\end{figure}

\begin{figure}[htbp]
  \centering
  \label{fig:a}\includegraphics{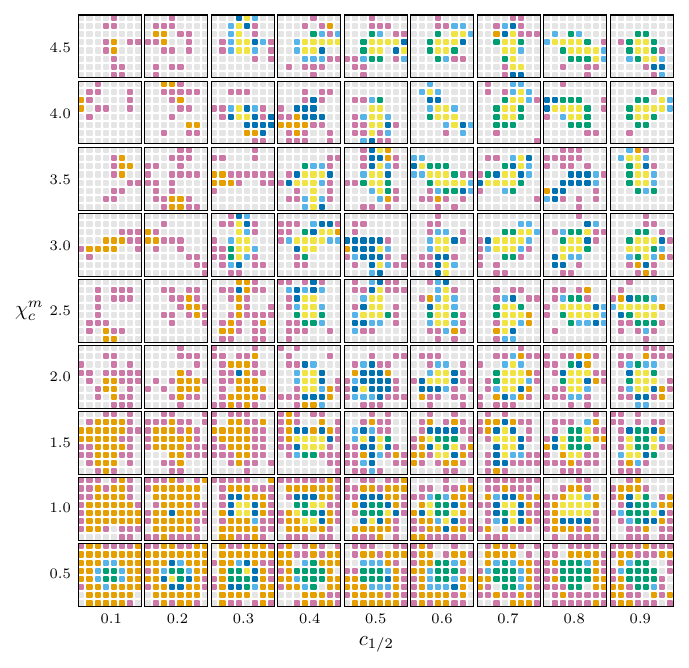}
  \caption{Categorisation of local neighbourhoods into 7 classes according to the local magnitudes of tumour cells, macrophages and necrotic cells.}
  \label{fig:local-magnitudes-7signatures}
\end{figure}

\begin{figure}[htbp]
  \centering
  \label{fig:a}\includegraphics{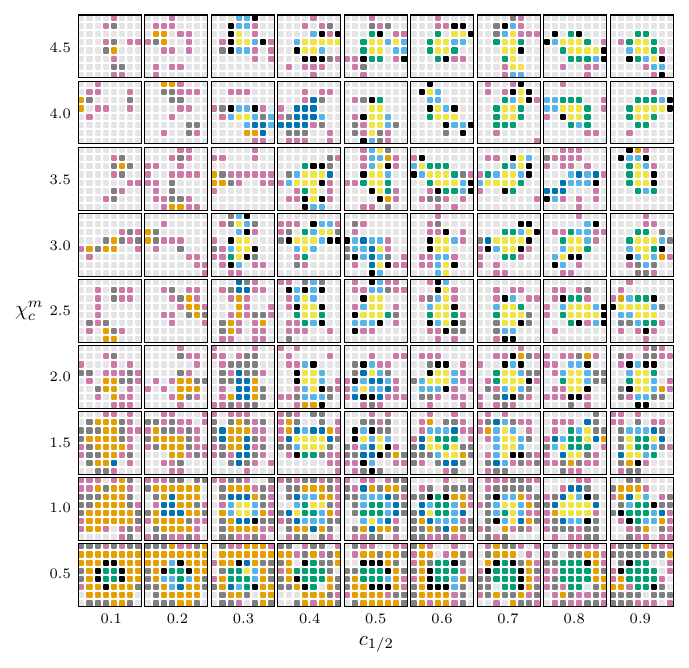}
  \caption{Categorisation of local neighbourhoods into 9 classes according to the local magnitudes of tumour cells, macrophages and necrotic cells.}
  \label{fig:local-magnitudes-9signatures}
\end{figure}

For simple magnitude feature vectors, the analysis suggests that $k \geq 5$ is needed to capture all significant cluster structure, and that for much larger $k$ any additional cluster structure may not be very significant. However, the increasing silhouette score indicates that cluster quality remains high. Hence, the analysis suggests that it could nonetheless be reasonable to increase $k$ further in order to obtain additional structure that may be more subtle. Results for $k=6$ are presented in Section~\ref{sec:local-magnitude} (see Figure~\ref{fig:local-across-schemes},~\ref{fig:signatures}), which yields some additional (radial) structure compared to $k=5$ (see Figure~\ref{fig:averages-magnitudes-5signatures},~\ref{fig:local-magnitudes-5signatures}) but remains simple enough for a subjective assessment.
For further reference, we present results for $k=7$ and $k=9$ in Figures~\ref{fig:local-magnitudes-7signatures} and~\ref{fig:local-magnitudes-9signatures}.

The analysis in Figure~\ref{fig:clustering-scores} further suggests that for simple cell count feature vectors, $k=5$ or $k=6$ is optimal, and for the more intricate cell count and distance vectors, $k=3$, $4$ or $7$ appear to be good candidates. In favour of direct comparabilty to magnitude features, Section~\ref{sec:local-magnitude} and Appendix~\ref{app:cellcounts} include results for $k=6$.

\section{Global Application: Previous Results}\label{app:global-previous}
\begin{figure}[htbp]
  \centering
  \label{fig:a}\includegraphics{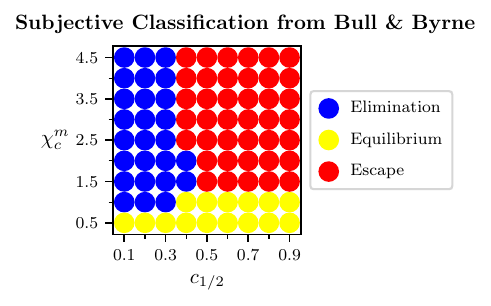}
  \caption{Subjective classification from \cite{Bull2023} of \(9 \times 9\) pairs of parameter values \(\chi_c^m\) and \(c_{1/2}\)).}
  \label{fig:results-bull}
\end{figure}

\begin{figure}[htbp]
  \centering
  \label{fig:a}\includegraphics{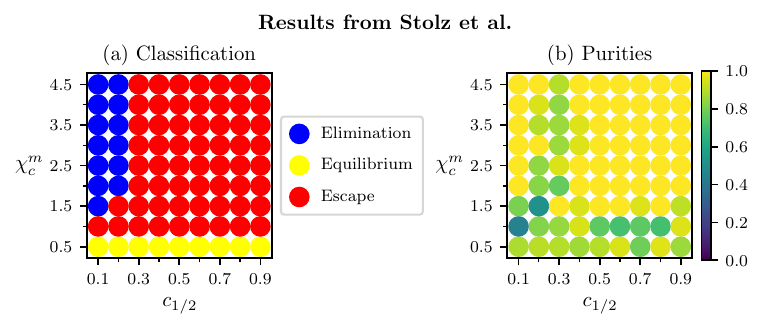}
  \caption{(a) Classification of \(9 \times 9\) pairs of parameter values \(\chi_c^m\)  and \(c_{1/2}\) from \cite{Stolz2024}. (b) Purity scores of the classification quantify how consistently simulation outcomes within one parameter scheme were assigned to the dominant tumour behaviour within the scheme. The plots were created using code from the \emph{MultiplexWitnessPH.py} file on the GitHub repository~\cite{githubMultiplexRelations}.}
  \label{fig:results-stolz}
\end{figure}

The work \cite{Bull2023} of Bull \& Byrne~2023 includes a subjective classification, here depicted in Figure~\ref{fig:results-bull}. Furthermore, the aforementioned classification results from Stolz \textit{et al.}\ 2024~\cite{Stolz2024} are depicted in Figure~\ref{fig:results-stolz}. Their work introduces a new Witness type complex on multispecies data which captures how one species relates to all the others. They construct this complex for multiple pairs of cell types on each of the $9 \times 9 \times 20$ simulation outcomes to extract topological features. The resulting feature vectors are clustered using $k$-means clustering for $k=3$. This yields a partition of all simulations into the three qualitative behaviours Elimination, Equilibrium and Escape. Assigning each parameter scheme to the behaviour most dominant among the 20 corresponding simulations gives the depicted classification of parameter schemes. Furthermore, they calculate the \emph{purity} of each parameter scheme, i.e., the number of outcomes assigned to the dominant behaviour relative to the total number of outcomes within the given parameter scheme. The authors of~\cite{Stolz2024} argue that these results cannot be reproduced using simple descriptor vectors consisting of cell counts and average distances of each cell type to the nearest blood vessel (for a comparison, see Figure 1 in SI of~\cite{Stolz2024}).

\subsection{Moran's I}
\begin{figure}[htbp]
  \centering
  \label{fig:a}\includegraphics{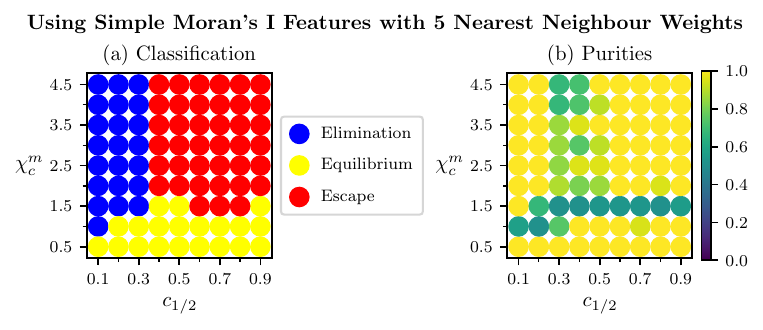}
  \caption{(a) Classification of \(9 \times 9\) pairs of parameter values \(\chi_c^m\)  and \(c_{1/2}\) using Moron's I features (see Equation~\ref{eq:morans-i} with $5$-nearest neighbour weights. (b) Purity scores of the classification quantify how consistently simulation outcomes within one parameter scheme were assigned to the dominant tumour behaviour within the scheme.}
  \label{fig:classification-moran-knn}
\end{figure}

\begin{figure}[htbp]
  \centering
  \label{fig:a}\includegraphics{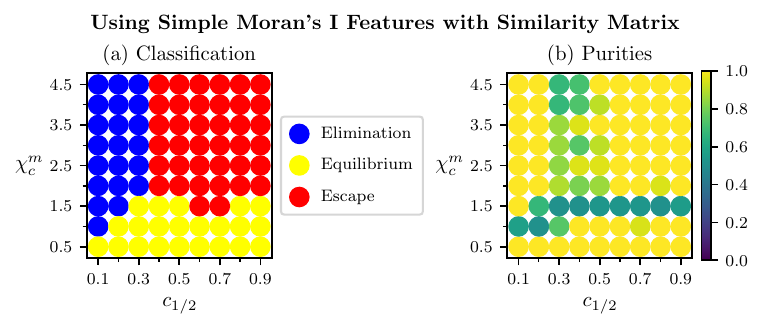}
  \caption{(a) Classification of \(9 \times 9\) pairs of parameter values \(\chi_c^m\)  and \(c_{1/2}\) using Moron's I features (see Equation~\ref{eq:morans-i} with weights obtained from similarity matrices. (b) Purity scores of the classification quantify how consistently simulation outcomes within one parameter scheme were assigned to the dominant tumour behaviour within the scheme.}
  \label{fig:classification-moran-similarity}
\end{figure}

\begin{figure}[htbp]
  \centering
  \label{fig:a}\includegraphics{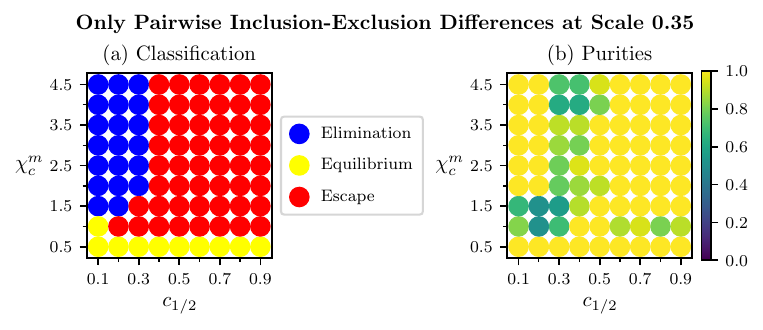}
  \caption{(a) Classification of \(9 \times 9\) pairs of parameter values \(\chi_c^m\)  and \(c_{1/2}\) using only pairwise inclusion-exclusion type differences of magnitudes (see Equation~\ref{eq:pairwise-inclusion-exclusion}). (b) Purity scores of the classification quantify how consistently simulation outcomes within one parameter scheme were assigned to the dominant tumour behaviour within the scheme.}
  \label{fig:classification-pairwise-differences}
\end{figure}

The inclusion-exclusion type differences of magnitudes used in the global application measure spatial correlation between multiple cell types. In particular, differences of the form $|X_1| + |X_2| - |X_1 \sqcup X_2|$ for a pair of species $X_1$ and $X_2$ encode their spatial relation. The same can be done using Moran's~I~\cite{moran_continous_1950}, a measure of spatial autocorrelation. While traditionally applied to spatial data with a continuous label, it can also be applied in this setting with a pair of species, by assigning the value 0 to one species and 1 to the other. Doing so for every pair of species in a simulation outcome yields an index $I(X_1,X_2) \in [-1,1]$ for each pair. Using the Moran's~I
\begin{equation}\label{eq:morans-i}
    (I(X_T,X_M),...,I(X_N,X_V))
\end{equation}
as feature vectors in our global pipeline yields a direct comparison to using just pairwise inclusion-exclusion type differences of magnitudes
\begin{equation}\label{eq:pairwise-inclusion-exclusion}
    (|X_T| + |X_M| - |X_T \sqcup X_M|,...,|X_N| + |X_V| - |X_N \sqcup X_V|).
\end{equation}
$I$ is set to 0 if one of $X_1$ and $X_2$ is empty. For the comparison we use two different sets of weights for Moran's I: $k$-nearest neighbour weights, as well as weights obtained from the similarity matrix (as intorduced in Section~\ref{sec:magnitude-definitions-properties}) using the full2W utility from the libpysal package~\cite{pysal_2007}. The former is a standard choice, while the latter is a natural choice in the context of magnitude.
The resulting classifications for the three sets of features are shown in Figures~\ref{fig:classification-moran-knn},~\ref{fig:classification-moran-similarity},~\ref{fig:classification-pairwise-differences}. Furthermore, the Moran's I features~\eqref{eq:morans-i} and magnitude features~\eqref{eq:pairwise-inclusion-exclusion} have a Pearson correlation greater than 0.99 for both considered choices of weights. These results confirm that the inclusion exclusion differences of magnitudes work exactly as intended. Magnitude inclusion-exclusion type differences recover spatial correlation between species, yet they appear to be picking up on slightly different patterns in boundary regions of the classifications compared to Moron's I. They further allow naturally for combination with singletons $|X_i|$ or inclusion-exclusion type differences with more than two species involved (see Equation~\eqref{eq:incl-excl}).

\newpage
\section{Figures}



\end{document}